\documentclass[a4paper]{amsart}

\usepackage{color}
\usepackage{graphicx}
\usepackage[normalem]{ulem}
\usepackage{morefloats}
\usepackage{hyperref}
\numberwithin{equation}{section}
\hypersetup{
    bookmarks=true,         
    unicode=false,          
    pdftoolbar=true,        
    pdfmenubar=true,        
    pdffitwindow=false,     
    pdfstartview={FitH},    
    pdftitle={My title},    
    pdfauthor={Author},     
    pdfsubject={Subject},   
    pdfcreator={Creator},   
    pdfproducer={Producer}, 
    pdfkeywords={keyword1, key2, key3}, 
    pdfnewwindow=true,      
    colorlinks=false,       
    linkcolor=green,          
    citecolor=green,        
    filecolor=green,      
    urlcolor=green,           
    urlbordercolor={1 1 1}  
}

\usepackage{amsmath, amsfonts, amsthm, amssymb}
\usepackage{fancyhdr}

\usepackage{biblatex}

\usepackage{mathrsfs}
\usepackage{blkarray}
\usepackage{amsthm}
\usepackage{mathtools}
\usepackage[title]{appendix}
\newtheorem{theorem}{Theorem}[section]
\newtheorem*{theorem*}{Theorem}

\newcommand{\Z}{\mathbb{Z}}

\theoremstyle{definition}

\theoremstyle{remark}

\begin{document}

\title[\resizebox{4.5in}{!}{\large  A NOTE ON A SPECTRAL CONSTANT ASSOCIATED  WITH AN ANNULUS}]{\large  A NOTE ON A SPECTRAL CONSTANT ASSOCIATED  WITH AN ANNULUS}

\author{ \small GEORGIOS TSIKALAS}
\address{DEPARTMENT OF MATHEMATICS AND STATISTICS, WASHINGTON UNIVERSITY IN ST. LOUIS, ST. LOUIS, MO, 63136}
\email{gtsikalas@wustl.edu} 
\keywords{Spectral set, $K$-spectral set, intersection of spectral sets. \\  \hspace*{0.34 cm} I thank Professor John E. M$^c$Carthy, my thesis advisor, for his invaluable advice during the \\ preparation of this work.}
\small
\begin{abstract}
    \small
 Fix $R>1$ and let $A_R=\{1/R\le |z|\le R \}$ be an annulus. Also, let $K(R)$ denote the smallest constant such that $A_R$ is a $K(R)$-spectral set for the bounded linear operator $T\in \mathcal{B}(H)$  whenever $||T||\le R$ and $||T^{-1}||\le R.$ We show that $K(R)\ge 2, \text{ for all } R>1. $ This improves on previous results by Badea, Beckermann and Crouzeix.
\end{abstract}
\maketitle

 \section{BACKGROUND}
 \large 
\par Let $X$ be a closed set in the complex plane and let $\mathcal{R}(X)$ denote the algebra of complex-valued bounded rational functions on $X,$ equipped with the supremum norm 
 $$||f||_{X}=\sup\{|f(x)|:x\in X \}.$$ Suppose that $T$ is a bounded linear operator acting on the (complex) Hilbert space $H$. Suppose also that the spectrum $\sigma (T)$  of $T$ is contained in the closed set $X$. Let $f=p/q\in\mathcal{R}(X)$. As the poles of the rational function $f$ are outside of $X,$ the operator $f(T)$ is naturally defined as $f(T)=p(T)q(T)^{-1}$ or, equivalently, by the Riesz-Dunford functional calculus (see e.g. \cite{Conway} for a treatment of this topic).
\par Recall that for a fixed constant $K>0,$ the set $X$ is said to be a \textit{$K$-spectral set} for $T$ if $\sigma (T)\subseteq X$ and the inequality $||f(T)||\le K||f||_X$ holds for every $f\in\mathcal{R}(X).$
The set $X$ is a \textit{spectral set} for $T$ if it is a $K$-spectral set with $K=1.$  Spectral sets were introduced and studied by von Neumann in \cite{vNeumann}, where he proved the celebrated result that an operator $T$ is a contraction if and only if the closed unit disk is a spectral set for $T$ (we refer the reader to the book \cite{Paulsen} and the survey \cite{badeaspectral} for more detailed presentations and more information on $K$-spectral sets).
\par We will be concerned with the case where $X=A_R:=\{1/R\le |z|\le R\}$ ($R>1$) is a closed annulus, the intersection of the two closed disks $D_1=\{|z|\le R \}$ and $D_2=\{|z|\ge 1/R \}$. Now, the intersection of two spectral sets is not necessarily a spectral set; counterexamples for the annulus were presented in \cite{Misra}, \cite{Towardatheory} and \cite{Minimalspectral}. However, the same question for $K$-spectral sets remains open (the counterexamples for spectral sets show that the same constant cannot be used for the intersection). Regarding the annulus in particular, Shields proved in \cite{Shields} that, given an invertible operator $T\in\mathcal{B}(H)$ with $||T||\le R$ and $||T^{-1}||\le R,$ $A_R$ is a $K$-spectral set for $T$ with $K=2+\sqrt{(R^2+1)/(R^2-1)}$. This bound is large if $R$ is close to $1$. In this context, Shields raised the question of finding the smallest constant $K=K(R)$ such that $A_R$ is $K(R)$-spectral. In particular, he asked whether this optimal constant $K(R)$ would remain bounded. \par 
This question was answered positively by Badea, Beckermann and Crouzeix  in \cite{IntersectionofSeveralDisks}, where they obtained that (for every $R>1$)
$$  \resizebox{0.99\hsize}{!}{$\frac{4}{3}<\gamma(R):=2(1-R^{-2})\prod_{n=1}^{\infty}\big(\frac{1-R^{-8n}}{1-R^{4-8n}}\big)^2\le K(R)\le 2+\frac{R+1}{\sqrt{R^2+R+1}}\le 2+\frac{2}{\sqrt{3}}.$}$$

\large
It should be noted that the quantity $\gamma(R)$ was numerically shown to be greater than or equal to $\pi/2$ (leading to the universal lower bound $\pi/2$ for $K(R)$) and it also tends to $2$ as $R$ tends to infinity.\par
Two subsequent improvements were made to the upper bound for $K(R)$: the first one in \cite{AnnulusasK-spectral} by Crouzeix and the most recent one in \cite{Crouzeix-SIAM} by Crouzeix and Greenbaum, where it was proved that $$ K(R)\le 1+\sqrt{2}, \hspace{0.3 cm}\forall R>1.$$
As for the lower bound, Badea obtained in \cite{K-spectralsetsbyBadea} the statement 
 $$\cfrac{3}{2}<2\frac{1+R^2+R}{1+R^2+2R}\le K(R), \hspace{0.3 cm}\forall R>1,$$
 where the quantity $2(1+R^2+R)/(1+R^2+2R)$ again tends to $2$ as $R$ tends to infinity. \\
 We improve the aforementioned estimates by showing that $2$ is actually a universal lower bound for $K(R)$:
\begin{theorem}\label{1.1}
Put $A_R=\{1/R\le |z|\le R \}$, for any $R>1$. Let $K(R)$ denote the smallest positive constant such that $A_R$ is a $K(R)$-spectral set for the bounded linear operator $T\in \mathcal{B}(H)$  whenever $||T||\le R$ and $||T^{-1}||\le R.$ Then, $$ K(R)\ge 2, \hspace{0.3 cm}\forall R>1.$$
\end{theorem} \vspace*{0.1 cm}
\small
\section{PROOF OF THEOREM \ref{1.1}}
\large
\begin{proof}
Fix $R>1$. For every $n\ge 2$, define 
$$g_n(z)=\frac{1}{R^n}\bigg(\frac{1}{z^n}+z^n\bigg)\in \mathcal{R}(A_R).$$
It is easy to see that 
\begin{equation}\label{1} ||g_n||_{A_R}=g_n(R)=1+\frac{1}{R^{2n}}. \end{equation}
To achieve the stated improvement, we will apply $g_n$ to a bilateral shift operator $S$ acting on a particular weighted sequence space $L^2(\beta)$.
First, define the sequence $\{\beta(k)\}_{k\in\Z}$ of positive numbers (weights) as follows:  
$$\beta (2ln+q)=R^q, \hspace{0.2 cm} \forall q\in\{0,1,\dots, n\},\forall l\in\Z; $$
$$\beta ((2l+1)n+q)=R^{n-q}, \hspace{0.2 cm} \forall q\in\{0,1,\dots, n\},\forall l\in\Z.$$ \\
Consider now the space of sequences $f=\{\hat{f}(k) \}_{k\in\Z}$ such that 
$$||f||^2_{\beta}:=\sum_{k\in\Z}|\hat{f}(k)|^2[\beta(k)]^2<\infty.$$
We shall use the notation  $f(z)=\sum_{k\in\Z} \hat{f}(k)z^k$ (formal Laurent series), whether or not the series converges for any (complex) values of $z$.
Our weighted sequence space will be denoted by 
$$L^2(\beta):=\{f=\{\hat{f}(k) \}_{k\in\Z}:||f||^2_{\beta}<\infty\}.$$
This is a Hilbert space with the inner product 
$$\langle f,g\rangle_{\beta}:=\sum_{k\in\Z} \hat{f}(k)\overline{\hat{g}(k)}[\beta(k)]^2. $$
Consider also the linear transformation (bilateral shift) $S$ of multiplication by $z$ on $L^2(\beta)$: 
                $$(Sf)(z)=\sum_{k\in\Z}\hat{f}(n)z^{n+1}.$$ 
                In other words, we have 
                $$\widehat{(Sf)}(n)=\hat{f}(n-1), \hspace{0.2 cm} \forall n\in\Z.$$
Observe that $$||S||=\sup_{k\in\Z}\frac{\beta(k+1)}{\beta (k)}=R $$
and $$||S^{-1}||=\sup_{k\in\Z}\frac{\beta(k)}{\beta (k+1)}=R.$$
Now, let $m\ge 3$ and define $h=\{\hat{h}(k) \}_{k\in\Z}\in L^2(\beta)$ by putting: 
$$\hat{h}(2ln)=\frac{1}{m}, \hspace{0.2 cm} \forall l\in\{0,1,2\dots,m^2\};$$
$$\hat{h}(k)=0, \hspace{0.2 cm} \text{in all other cases}.$$
We calculate 
$$||h||^2_{\beta}=\sum_{l=0}^{m^2}\frac{1}{m^2}[\beta(2ln)]^2=\sum_{l=0}^{m^2}\frac{1}{m^2}\cdot 1^2=\frac{m^2+1}{m^2},$$
hence
\begin{equation}\label{2}  ||h||_{\beta}=\frac{\sqrt{m^2+1}}{m}.\end{equation}
Also, put $f=(S^{-n}+S^n)h$ and notice that 
$$||(S^{-n}+S^n)h||^2_{\beta}=||f||^2_{\beta} $$ $$\ge \sum_{l=1}^{m^2}|\hat{f}((2l-1)n)|^2[\beta((2l-1)n)]^2$$ $$=\sum_{l=1}^{m^2}\bigg(\frac{2}{m}\bigg)^2 R^{2n}$$ $$=4R^{2n}.$$
Thus,
\begin{equation}\label{3}  ||(S^{-n}+S^n)h||_{\beta}\ge 2R^n . \end{equation}
Using (\ref{1}), (\ref{2}) and (\ref{3}), we can now write 
 $$K(R)\ge \frac{||g_n(S)||}{||g_n||_{A_R}} $$ $$=\cfrac{1}{R^n}\cdot\cfrac{||S^{-n}+S^n||}{1+R^{-2n}}$$ $$\ge \cfrac{1}{R^n+R^{-n}}\cdot\cfrac{||(S^{-n}+S^n)h||_{\beta}}{||h||_{\beta}}$$ $$\ge \cfrac{1}{R^n+R^{-n}}\cdot\cfrac{2R^n}{\frac{\sqrt{m^2+1}}{m}}.$$
Letting $m\to\infty,$ we obtain 
$$K(R)\ge \cfrac{1}{R^n+R^{-n}}\cdot\cfrac{2R^n}{1}=\cfrac{2R^n}{R^n+R^{-n}}\xrightarrow{n\to\infty}2,\hspace{0.2 cm} \text{ as } R>1.$$
The proof is complete. 

\end{proof}

\printbibliography

\end{document}